\documentclass{article}
\usepackage[utf8]{inputenc}

\usepackage{graphicx}
\usepackage{amsmath}
\usepackage{caption}
\usepackage{verbatim}
\usepackage{indentfirst}
\usepackage{titlesec}
\titlelabel{\thetitle.\quad}

\begin{document} 

\title{New Method for Building Geodesic Lines in Riemann Geometry}
\author{Victor I. Pogorelov, Ph.D.}
\date{}

\maketitle

\noindent New perspective form of equations for geodesic lines in Riemann Geometry was found. 


\section{Introduction}
This work is a result of the development of the Method of Differential Coordinates (MDC) that has been designed as a method for solving differential equations \cite{A} \cite{B}. As we know, the majority of physical processes are described by differential equations in partial derivatives of second order linear with respect to second derivatives. The analytical stage of solving problems with equations of this type was important until recently when computers gained capacity to use these equations without preliminary analytical processing. But even nowadays, direct numerical (grid) calculation is not always possible and is not always expedient, especially because it masks the physics of the process and does not explain why the result of the calculation is what it is. That is why analytical methods (and MDC is partially one of them) are still very important. However, the present work is not devoted to the subject of differential equations, although its mathematical foundation is laid there. Previous tasks required the use of attributes of variational computation (geodesic lines), which is understandable, even familiar, in a sense. However, the specifics of MDC allow for a new outlook at what seems to be well-established stereotypes.

In the process of using MDC, the computer simultaneously builds couples of geodesic lines that are connected to each other by the very process of their building. The method of their building does not use any equations that used to be employed for computation of separate geodesic lines. \textbf{The transition from conjugate calculations of several close geodesic lines to calculation of one solitary line is, in fact, the task of the present work. It is achieved by extreme convergence, integration of a couple of geodesic lines.} Let us note that the possibility of determining extremals’ trajectories by combined analysis of course of their pairs has been shown by Christiaan Huygens explanation of light wave flow forms long before creation of Riemannian geometry \cite{C}.

But the main idea of this work – the one that is also the basis for MDC – is connected with the introduction and application of differential coordinates.

\section{Differential Coordinates and their Properties }
The fundamentals of MDC are described in a number of the author’s works, \cite{A} in particular. Here are just the summaries, to avoid repetitions and complication of the text.

MDC grounds in the fact that a derivative, as a concept, needs only availability of the argument’s differential but not necessarily its existence as a non-zero value. Therefore, if in a coordinate space $x$ (that is $x_1, x_2, \dots x_n$) is given a differential field $dy$ (that is $dy_1 , dy_2 , \dots dy_n$), where
\begin{center}
\begin{equation}
    dy_i = \sum \alpha_{ij}(x)dx_j,
\end{equation}
\end{center}

\noindent then a derivative of any function $f(x)$ for any $dy_i$ can be determined as

\begin{center}
\begin{equation}
    \frac{Df(x)}{Dy_i} = \lim_{x' \to x} \frac{f(x)-f(x')}{\int_{x'}^{x} dy_i}
\end{equation}
\end{center}

\noindent without requirement that differentials $dy_i$ are complete, i.e. without requirement of existence of values $y_i$ (the only important condition is equality of all $dy_j =0$ at $j\neq i$ on the path of integration from $x'$ to $x$ along coordinates’ line of coordinate $dy_i$). Taking this factor into account expands the concept of the derivative, and, at that, operator $\frac{\partial}{\partial x_j}$ becomes a special case of the operator $\frac{D}{Dy_j}$. Summation in (1) and below is conducted by repeating indexes from 1 to $n$. 

Squaring the elementary distance in coordinates $dy_i$, that is $\sum dy_k^2 =  ds^2$, after applying here (1) will look like:

\begin{equation}
    ds^2 = \sum \rho_{ij}(x)dx_i dx_j,
\end{equation}

\noindent where $\rho_{ij} =\sum \alpha{k_i}\alpha{k_j}$. It is evident that this way the original formula of Riemann geometry is created. If (3) is given initially, then $dy_i$ through $\rho_{ij}$ is found by applying in (3) the expression of the reverse transition - from $dx_j$ to $dy_i$. 

\begin{equation}
    dx_i = \sum \alpha_{ik}^* dy_k,
\end{equation}

\noindent and then $\alpha_{ik}^*$ follow from the system of equations $\sum \rho_{ij}\alpha_{ik}^*\alpha_{j\nu}^* = \delta_{k\nu}$.

Coordinates $dy_i$ can be converted to $dy_j'$ (which is isometric to $dy_i$) by conventional way of rotation: $dy_i'  = \sum \beta_{ij} dy_j$, then from condition $\sum dy_i^2 =  \sum dy_j'^2$ follows $\sum \beta_{ij} \beta_{ik} =  \delta_{jk}$. It is easy to see that operator $\Delta y$ is covariant in such transformations, it retains its form.

It is shown in \cite{A} that in every point $x$ for every function that has any non-zero first derivatives in this point, is possible to define the gradient on the basis of $dy_j'$ as the only coordinate $dy_j'$ along which this function has a non-zero derivative. Note that in the case of complex functions the additional condition is that the function is analytic.

Also, we can show from (2) that the gradient’s square, that is $\sum (\frac{Df}{Dy_i})^2$, is an invariant for rotation transformations, along with $ds^2$.

The important property of the differential coordinates is that if in every point of solution they are rotated in such a way so one of them would be directed along function’s gradient, then the differential equation in partial derivatives for this function, written via given coordinates, becomes an equation similar to an ordinary (one-dimensional) differential equation.

\section{Geodesic Lines and their Properties}
If every element of the path in a given space is correlated with length $ds$, then lengths of all possible curves (including extreme, geodesic curves) are implicitly defined. In papers \cite{A} \cite{B}, analysis of such curves was conducted for solving differential equations. Here we have a different task, but both topics are equally dealing with geodesic curves’ properties.

Let’s assume that in $n$-dimensional space with specified parameters of the medium, i.e. with a specifically defined equation (3), from point $x_0$ comes a system of geodesic lines, and thus $x$ is given a function in every point \newline \noindent $S(x_0 , x) =  \int_{x_0}^x ds$ , where  integration is from point $x_0$ - the source of  geodesic lines – to the given point $x$ along the corresponding geodesic curve. Then closed $(n-1)$-dimensional hypersurfaces surrounding point $x_0$ where $S$=const will also be defined. Additionally, the system of coordinates $dy_i'$ will be partially defined. If, for example, $dy_1'$ is oriented along gradient of function $S(x)$, then  all other $dy_i'$ will fall into the area $S=const$.

Let us note two properties that occur in the system of geodesic lines in question, considering, for example, a two-dimensional case, when hypersurfaces $S=const$ are one-dimensional curves (Fig.1).

Increase of value $S$ (which is dictated by integration of any two fixed hypersurfaces $S=const$) will be the same along all geodesic lines. After all, this increase, i.e. the difference of integrals $\int_{x_0}^x ds$, for example between lines $S_3 = const$  and  $S_4 = const$, will be $S_4 - S_3$, no matter which geodesic line we take. That is why given differences of integrals on segments between two fixed hypersurfaces $S=const$ are the same for all geodesic lines.

Furthermore, the gradient of value $S(x)$ in “$y$-metrics” will be always directed along geodesic line. Indeed, the neighborhood area of current point $x$ without its one-dimensional subspace – vector $dS$, where $S=const$, – has $(n-1)$ dimensions. And the neighborhood area, that has dimension along direction of gradient $S$ extracted and is left with only part of space where $S=const$, is also $(n-1)$-dimensional; i.e. these two $(n-1)$-dimensional subspaces coincide. Therefore, two one-dimensional subspaces – vector $dS$ and gradient $S$ – have to coincide too, which means coincidence of directions of gradient $S(x)$ and the geodesic line.

Fig. 1 shows how geodesic lines curve in the area between $S_1$  and  $S_2$  assuming that functions $S(x_0, x)$ and $\alpha(x)$  are either real or are projections of real parts of complex $S$, $\alpha$ onto space of real values $x_i$. There two close to each other geodesic lines cross also close to each other hypersurfaces $S_1$ and $S_2$. Let’s define the geodesic lengths that are count off point $x_0$ in “$x$-metrics” as $r_1 x$, $r_2 x$, and the angle between axis $x_2$ and $d r_i$ as $F_i$. If we turn coordinates $x$ in some point of the geodesic line in such a way that the new differential $dx_1'$ is directed along this line, then we get a corresponding factor $\alpha(x)$ for equation that connects the element $dy_1'= dS$ along the geodesic line with $(dx_1')_i = dr_i$ , and equation (1) for $dy_1'$ will become $dS = \alpha(x) dr_i$. Increase of value $S$ between $S_1$=const and $S_2=$const for both curves is $\delta S = S_2 - S_1$, i.e. values $\delta S = \alpha(x')\delta r_1$ and $\delta S= \alpha (x'') \delta r_2$ for both of them are the same; but as $\alpha(x')$  and $\alpha(x'')$ are different, the values $\delta r_1$ and $\delta r_2$ are different too.

Without knowing how the angle $F$ of vector $dr$ of every geodesic will change between $S_1$ and $S_2$, we can draw zones’ borders on their ends in $S_1$ that will correspond with given value $\delta S$.

For simplicity, on Fig.1 these borders are shown as circles, although it is possible for them to take shape of deformed circles (like ellipses) if $\alpha$ will depend also on direction of vector $dr$. But that has no basic importance. The curve $S_2$ (i.e. hypersurface $S_2$) can only touch these zones’ borders or there would be shorter ways to $S_2$ than geodesic lines, and that is impossible. Based on that, we can calculate $\delta F$ – variation of angular orientation of elements $dr_i$ in the ends of given couple of geodesics \cite{A} after passing the distance from $S_1$ to $S_2$.

Small quantity $\delta F$ is generated by two small quantities: \newline \noindent $\lambda =[(x_2''-x_2')^2+(x_1''-x_1')^2]^{\frac{1}{2}}$ (the result of initial divergence of geodesic lines in point $x_0$) and $\delta r_1$ (the result of divergence of $S_2$ from $S_1$). Up to second order values, small quantity $\delta F$ is connected to $\lambda$ and $\delta r_1$ with equation \cite{A} \cite{B}

\begin{center}
\begin{equation*}
    \sin(\delta F) = \frac{\delta r_1 \{ \frac{\alpha (x')}{\alpha (x'')}-1\} }{\lambda}.
\end{equation*}
\end{center}

\noindent Expanding $sin(\delta F)$ in a row and tending $\lambda$ and $\delta r_1$ by turns to zero, we have:

\begin{center}
\begin{equation}
    \frac{dF}{dr} = \frac{\alpha '}{\alpha },
\end{equation}
\end{center}

\noindent where $\alpha'$ is a derivative of $\alpha$ with respect to orthogonal to $\delta r_1$ in “$x$-metrics” in direction of gradient $\alpha$ (the result of the operation $\lambda \to 0$). From here, knowing that $\frac{dx_2}{dr}= \cos F$ , differentiating it one more time in $r$, and taking in account (5), we find:

\begin{center}
\begin{equation}
    \frac{d^2x_2}{dr^2} = - \frac{\alpha '}{\alpha} \sqrt{1-\Big( \frac{dx_2}{dr} \Big)^2}.
\end{equation}
\end{center}

Naturally, the analog of this formula may be written also for $\frac{d^2x_1}{dr^2}$. In principle, the root in the right part of (6) can be substituted for $\frac{dx_1}{dr}$, and then (6) will become linear with respect to derivative of $x_1$. But the value of the variant (6) lies in dealing only with derivatives in one coordinate, the very coordinate that interests us at this moment.

When $n>2$, a two-dimensional space should be allocated with respect to vectors $dr$ and $dx_j$ in the neighborhood of the geodesic line. Then temporarily introducing second (intermediary) coordinate of type $x_t$ , we can obtain an equation (6) for $x_j$. As a result, the geodesic will be defined by $n$ equations

\begin{center}
\begin{equation}
    \frac{d^2x_j}{dr^2} = - \frac{\alpha '}{\alpha} \sqrt{1- \Big( \frac{dx_j}{dr} \Big)^2},
\end{equation}
\end{center}

\noindent where the derivative from $\alpha$ is always taken in a crosswise direction to $dr$ in an element of a two-dimensional surface that was created by vectors $dr$ and $dx_j$. Conceptually, the gradient $\alpha$ curves the geodesic in the plane of the gradient’s “pressure” on it.

Formula (7) has a very clear meaning. The square root is the cosine of an angle between the gradient of function $\alpha (x)$ and vector $dr$ normal. It is evident that gradient’s component along $dr$ does not change its alignment, and orthogonal to $dr$ component pulls this vector to itself.

Value $\alpha (x)$ that appears in (6) and (7) is introduced by connecting $dy_1'$ and $dr$ with the equation $dy_1' = \alpha (x) dr$. But $ds^2 = (dy_1')^2 = \alpha^2 dr^2$ . That is why if we are given a particular orientation of vector $dr$ (and particular values $\frac{dx_j}{dr}$ along with it) and also space metrics are predetermined by formula (3), then by placing $dx_k=\frac{dx_k}{dr}dr$ in it and taking into account the equation $ds^2 = \alpha^2dr^2$, we have

\begin{center}
\begin{equation*}
    \alpha = \sqrt{\sum \rho_{ij} (x) (\frac{dx_i}{dr})(\frac{dx_j}{dr})}.
\end{equation*}
\end{center}

Clearly, the absence of additional cumbersome calculations makes (7) very convenient to use. Now, let’s remember that Fig. 1 was built for real terms $S$ and $\alpha$. But if equation (6) and function $S(x)$ are analytical, than generalization (7) for complex $\alpha$, $x_i$ have to fully satisfy conditions of extremum (minimum) when its real term satisfies these conditions. At that, functions and their derivations maintain their values, and that secures the equations formed from them. As is well known, Riemannian geometry analogue of (7) which is grounded in customary conceptions of parallel transport of a vector along geodesic line, looks like

\begin{center}
\begin{equation}
    \frac{d^2x_i}{ds^2} + \sum \Gamma_{\alpha \beta}^i (\frac{dx_\alpha}{ds}) (\frac{dx_\beta}{ds}) = 0,
\end{equation}
\end{center}

\noindent where $\Gamma_{\alpha \beta}^i$ is the Cristoffel symbol which is denoted in cumbersome manner by $\rho_{ij}$ and their derivatives with respect to $x_v$ . Nowadays, (8) has become one of the pillars for the theory of vector and tensor analysis. The conclusion (7), drawn from a different principle, is simpler, and the formula itself is more understandable. 

What comes to mind when we compare these formulas? They are interchangeable, but (7) might have an advantage over (8), as it can be more convenient to use. In particular, it can be used to illustrate and explain results of application of (8). Let’s note that (7) and (8) were tested during deduction of Snell’s law, which is well known in radio-physics. Both of them have proved the law to be true. (7) has proved the law immediately, and (8) after long laborious computations. It is evident that (7) is highly promising for theoretical and practical physics. Moreover, because of its clarity, it is very convenient for fast qualitative valuations prior to final precise calculations.

Let's emphasize the following: (8) had served as a basis for the creation of general theory of relativity; perhaps (7) will prove to be helpful for development in this direction.

\newpage

\begin{figure}[h]
    \centering
    \includegraphics[scale = .75]{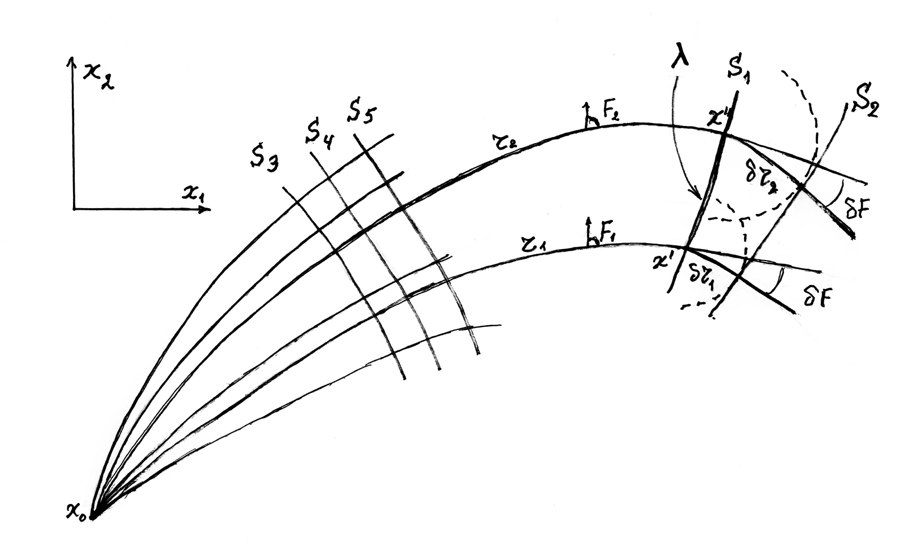}
    \caption{}
\end{figure}

\bibliographystyle{amsplain}

\end{document}

\title{Новый метод построения геодезических линий в римановой геометрии.}
\author{В. И. Погорелов}
\date{}

\selectlanguage{russian}

\maketitle

\selectlanguage{english}

\noindent New perspective form of equations for geodesic lines in Riemann Geometry was found. 

\noindent {\bf Keywords:} Riemann Geometry, geodesic lines, differential coordinates
\newline

\noindent\rule[0.5ex]{\linewidth}{.05pt}
\selectlanguage{russian}

\noindent Найдена новая форма уравнений для геодезических линий римановой геометрии, перспективная для их изучения. 

\noindent {\bf Ключевые слова:} геометрия Римана, геодезические линии, дифференциальные координаты

\section{Введение}
Данная работа явилась результатом развития метода дифференциальных координат (МДК), который  возник как способ решения дифференциальных уравнений \cite{A} \cite{B}. Как известно, большинство физических процессов описываются дифференциальными уравнениями в частных производных второго порядка, линейных относительно старших (2–х) производных. Аналитический этап решения задач с такими уравнениями до недавнего времени был очень актуален, когда техника вычислительных машин не позволяла их применять сразу без предварительной аналитической обработки задачи. Но и сейчас прямой численный (сеточный) расчёт не всегда возможен и не всегда целесообразен, тем более, что он маскирует физику процесса, не объясняет, почему результат счёта таков, какой он есть. Поэтому аналитические методы (к которым частично относится и МДК) и поныне играют большую роль. Впрочем, настоящая работа посвящена не тематике дифференциальных уравнений, хотя её математический фундамент заложен именно там. Прежние задачи привели к необходимости использования атрибутов вариационного исчисления – к геодезическим линиям, что в принципе понятно, даже в какой-то степени привычно. Однако, специфика МДК позволила по новому взглянуть на казалось бы уже прочно сложившиеся стереотипы. 
В ходе применения МДК машина строит сразу пары  геодезических линий, сопряжённых друг с другом самим процессом построения. Причём способ этого построения не использует уравнений, ранее применявшихся при расчётах отдельных геодезических. \textbf{Процесс перехода от сопряжённых расчётов нескольких близких геодезических линий к расчёту одной обособленной линии является задачей, которой фактически и посвящена настоящая работа. Это достигается путём предельного сближения, слияния пары геодезических.}  Следует отметить, что возможность определения траекторий экстремалей путём совместного анализа хода их пар показана ещё Х. Гюйгенсом при объяснении форм пути светового луча задолго до времени создания римановой геометрии \cite{C}. 

Но основная  идея настоящей работы, положившая и начало МДК, была связана с вводом и использованием дифференциальных координат.

\section{Дифференциальные координаты и их свойства}
Основы МДК изложены в ряде работ автора, в частности в \cite{A}. 
Здесь приводятся по возможности лишь их конечные выводы чтобы избежать повторений и не усложнять данный текст.

Метод МДК опирается на тот факт, что для  производной как понятия нужно лишь наличие  дифференциала аргумента, а существование его самого в виде  некой отличной от нуля величины не является необходимым. Поэтому если в пространстве координат $x$ (т.е. $x_1, x_2, \dots x_n$) задано дифференциальное поле $dy$ (т.е. $dy_1 , dy_2 , \dots dy_n$), где
\begin{center}
\begin{equation}
    dy_i = \sum \alpha_{ij}(x)dx_j,
\end{equation}
\end{center}

\noindent то можно определить производную любой функции $f(x)$ по любому  $dy_i$ как

\begin{center}
\begin{equation}
    \frac{Df(x)}{Dy_i} = \lim_{x' \to x} \frac{f(x)-f(x')}{\int_{x'}^{x} dy_i}
\end{equation}
\end{center}

\noindent не требуя того, чтобы дифференциалы $dy_i$ были полными, т.е. не требуя существования самих величин $y_i$ (важно лишь условие равенства всех $dy_j =0$ при $j\neq i$ на пути интегрирования от $x'$ до $x$ вдоль координатной линии координаты $dy_i$). Учёт этого факта расширяет понятие производной, а оператор $\frac{\partial}{\partial x_j}$ становится частным случаем оператора $\frac{D}{Dy_j}$. (Суммирование в (1) и ниже ведётся по повторяющимся индексам от 1 до $n$).  

Квадрат элементарного расстояния в координатах $dy_i$, т.е. $\sum dy_k^2 =  ds^2$, после подстановки сюда (1)  будет выглядеть как:

\begin{center}
\begin{equation}
    ds^2 = \sum \rho_{ij}(x)dx_i dx_j,
\end{equation}
\end{center}

\noindent где $\rho_{ij} =\sum \alpha{k_i}\alpha{k_j}$. Видно, что таким образом создаётся исходная формула римановой геометрии. Если же (3) задано изначально, то $dy_i$ через $\rho_{ij}$ находится подстановкой в (3) выражения обратного перехода – от $dx_j$ к $dy_i$. 

\begin{center}
\begin{equation}
    dx_i = \sum \alpha_{ik}^* dy_k,
\end{equation}
\end{center}

\noindent далее $\alpha_{ik}^*$ следуют из системы уравнений $\sum \rho_{ij}\alpha_{ik}^*\alpha_{j\nu}^* = \delta_{k\nu}$.

Координаты $dy_i$ могут преобразовываться к изометричным им $dy_j'$ обычным способом поворота: $dy_i'  = \sum \beta_{ij} dy_j$, далее из условия $\sum dy_i^2 =  \sum dy_j'^2$ вытекает $\sum \beta_{ij} \beta_{ik} =  \delta_{jk}$. Легко видеть, что при подобных преобразованиях оператор $\Delta y$ является ковариантным, сохраняет свой вид.

В \cite{A} показано, что в любой точке $x$ у любой функции, имеющей в этой точке какие-либо отличные от нуля первые производные, можно на основе $dy_j'$ определить градиент как единственную из координат   $dy_j'$ вдоль которой у этой функции существует не нулевая производная.  Заметим, что если речь идёт о комплексных функциях, то дополнительным условием является то, чтобы функция была аналитической.

Из (2) можно также показать, что квадрат градиента, т.е. $\sum (\frac{Df}{Dy_i})^2$, для преобразований поворота является инвариантом, как и $ds^2$.

Важным свойством дифференциальных координат является то, что если они в каждой точке области решения повёрнуты так, чтобы одна из них была направлена вдоль градиента функции, то дифференциальное уравнение в частных производных для этой функции, записанное через данные координаты, превращается в уравнение, близкое к обыкновенному (одномерному) дифференциальному уравнению.

\section{Геодезические линии и их свойства}
 Если каждому элементу пути в рассматриваемой области пространства сопоставлена длина  $ds$, то таким образом неявно определены длины всех возможных кривых, в том числе и экстремальных, геодезических. В работах \cite{A} \cite{B}, анализ таких кривых проводился для решения дифференциальных уравнений. Здесь  ставится иная задача, но обе тематики одинаково затрагивают свойства геодезических кривых.

Допустим, что в $n$-мерном пространстве с заданными параметрами среды, т.е. с конкретно определённым уравнением (3) , из точки $x_0$ исходит система геодезических линий, и таким образом в каждой точке $x$ задана функция  $S(x_0 , x) =  \int_{x_0}^x ds$ ,где интегрирование ведётся от т. $x_0$ – источника геодезических до данной точки $x$ вдоль соответствующей геодезической кривой. Тогда будут определены и  окружающие т. $x_0$ замкнутые $(n-1)$-мерные гиперповерхности, на которых $S=const$. Также частично определится и система координат  $dy_i'$. Если мы положим, что  $dy_1'$ ориентирована вдоль градиента функции $S(x)$, то все остальные  $dy_i'$ попадут в область $S=const$.

Отметим два свойства, возникающих у рассматриваемой нами системы геодезических линий, используя для иллюстрации двумерный случай, когда гиперповерхностями  $S = const$ будут одномерные кривые (Рис. 1).

Вдоль всех геодезических набег величины $S$ между двумя любыми фиксированными кривыми $S=const$ будет одинаков. Ведь он, т.е. разность интегралов $\int_{x_0}^x ds$, например между линиями $S_3 = const$ и $S_4 = const$, равен $S_4 - S_3$, о какой бы из геодезических линий не шла речь.   Поэтому данные разности интегралов на отрезках между двумя фиксированными гиперповерхностями $S=const$ одни и те же для всех геодезических линий.

Кроме того, градиент величины $S(x)$ в  « y-метрике» будет направлен всегда вдоль геодезической линии. Действительно, зона окрестности текущей точки $x$ без своего одномерного подпространства – вектора $dS$, на которой $S=const$, имеет $(n-1)$ измерений. И зона  окрестности, у которой удалено измерение вдоль направления градиента $S$ и осталась лишь та часть пространства, где $S=const$,тоже $(n-1)$-го измерения. Т. е. эти два $(n-1)$-мерных подпространства совпадают. Следовательно, должны совпасть и два одномерных подпространства – вектор $dS$ и градиент $S$, что и означает совпадение направлений градиента $S(x)$ и геодезической.

Рис. 1 иллюстрирует механизм искривления геодезических линий на участке между $S_1$ и $S_2$  в предположении, что функции $S(x_0, x)$ и $\alpha(x)$ или реальны, или являются проекциями реальных частей  комплексных $S$, $\alpha$ на пространство реальных значений $x_i$. Там две близкие геодезические линии пересекают  тоже близкие между собою гиперповерхности $S_1$ и $S_2$. Длины геодезических, отсчитываемые от т. $x_0$ в «$x$ -метрике», обозначим как  $r_1 x$, $r_2 x$, а угол между осью  $x_2$ и $d r_i$ - как $F_i$. Если координаты $x$  в  какой-либо точке геодезической повернуть так, чтобы новый дифференциал $dx_1'$ был направлен вдоль этой линии, то возникнет  соответствующий множитель $\alpha(x)$ у уравнения, связывающего элемент $dy_1'= dS$ вдоль геодезической с $(dx_1')_i = dr_i$, а уравнение (1) для $dy_1'$ приобретёт вид $dS = \alpha(x) dr_i$. Набег величины $S$  между $S_1=const$ и $S_2=const$ для обеих кривых равен $\delta S = S_2 - S_1$. Т.е. у них обеих значения $\delta S = \alpha(x')\delta r_1$ и $\delta S= \alpha (x'') \delta r_2$ одинаковы, но поскольку  $\alpha(x')$  и $\alpha(x'')$ различны, то различны также величины $\delta r_1$ и $\delta r_2$.

Не зная заранее, как между $S_1$ и $S_2$ изменится угол $F$ у вектора $dr$ каждой из геодезических, мы можем на концах их на $S_1$ начертить границы зон, соответствующих заданной величине $\delta S$.

На рисунке для простоты эти границы взяты в форме кругов, хотя и не исключено, что они могут принимать вид деформированных окружностей (вроде эллипсов), если $\alpha$ будут зависеть и от направления вектора $dr$. Но это не имеет принципиального значения. Кривая $S_2$ (т.е. гиперповерхность  $S_2$) может только касаться границ этих зон, иначе бы возникли и  более короткие пути к $S_2$ чем геодезические, что исключено. На  этой основе удаётся сделать расчёт $\delta F$ – изменения угловой ориентации элементов $dr_i$ на концах данной пары геодезических \cite{A} в результате прохождения пути от $S_1$ до $S_2$.

Малая величина $\delta F$ порождается двумя малыми величинами: \newline \noindent $\lambda =[(x_2''-x_2')^2+(x_1''-x_1')^2]^{\frac{1}{2}}$ (результат первичного расхождения геодезических в т. $x_0$) и $\delta r_1$ (результат отклонения $S_2$ от $S_1$). С точностью до величин второго порядка малости $\delta F$ связана с $\lambda$ и $\delta r_1$ уравнением \cite{A} \cite{B}

\begin{center}
\begin{equation*}
    \sin(\delta F) = \frac{\delta r_1 \{ \frac{\alpha (x')}{\alpha (x'')}-1\} }{\lambda}.
\end{equation*}
\end{center}

\noindent Разлагая $sin(\delta F)$ в ряд и устремляя поочерёдно $\lambda$ и $\delta r_1$ к нулю, имеем

\begin{center}
\begin{equation}
    \frac{dF}{dr} = \frac{\alpha '}{\alpha },
\end{equation}
\end{center}

\noindent где $\alpha'$ -   производная от $\alpha$ по ортогонали к $\delta r_1$ в «$x$ -метрике» в направлении градиента  $\alpha$ (результат операции $\lambda \to 0$). Отсюда, зная, что $\frac{dx_2}{dr}= \cos F$ , дифференцируя это ещё раз по $r$ и учитывая (5), находим

\begin{center}
\begin{equation}
    \frac{d^2x_2}{dr^2} = - \frac{\alpha '}{\alpha} \sqrt{1-\Big( \frac{dx_2}{dr} \Big)^2}.
\end{equation}
\end{center}

Естественно, аналог этой формулы можно написать и для $\frac{d^2x_1}{dr^2}$. В принципе, корень в правой части (6) можно заменить на $\frac{dx_1}{dr}$, и тогда (6)  станет линейной относительно производной от $x_1$. Но  ценность использованного варианта (6) заключается в том, что там фигурируют только производные от одной координаты, причём той, которой мы в данный момент интересуемся.

Когда $n>2$, в окрестности геодезической следует выделить двумерную область, задаваемую векторами $dr$ и $dx_j$. Далее в ней, вводя временно вторую (вспомогательную) координату типа $x_t$, можно так же  для  $x_j$ получить уравнение (6). В итоге геодезическая будет задаваться $n$ уравнениями    

\begin{center}
\begin{equation}
    \frac{d^2x_j}{dr^2} = - \frac{\alpha '}{\alpha} \sqrt{1- \Big( \frac{dx_j}{dr} \Big)^2},
\end{equation}
\end{center}

\noindent где производная от $\alpha$ берётся поперёк к $dr$ в элементе двумерной поверхности, созданной векторами $dr$ и $dx_j$. Образно говоря, градиент $\alpha$ искривляет геодезическую в плоскости своего «давления» на неё. 

Формула (7) имеет очень прозрачный смысл. Корень - это косинус угла между градиентом функции $\alpha (x)$ и нормалью к вектору $dr$. Видно, что компонента градиента вдоль $dr$ не изменяет его ориентацию, а компонента ортогональная $dr$ притягивает этот вектор к себе.

В (6) и (7) фигурирует величина $\alpha (x)$ , которая вводится уравнением связи $dy_1'$ и $dr$ как $dy_1' = \alpha (x) dr$. Однако, $ds^2 = (dy_1')^2 = \alpha^2 dr^2$. Поэтому, если задана конкретная ориентация вектора $dr$ (а с ним и конкретные значения величин $\frac{dx_j}{dr}$) и также задана метрика пространства формулой (3), то, подставляя в неё $dx_k=\frac{dx_k}{dr}dr$ и  учитывая уравнениеto $ds^2 = \alpha^2dr^2$, имеем

\begin{center}
\begin{equation*}
    \alpha = \sqrt{\sum \rho_{ij} (x) (\frac{dx_i}{dr})(\frac{dx_j}{dr})}.
\end{equation*}
\end{center}

Понятно, что использовать (7) очень удобно, так как это не требует дополнительных сложных вычислений. Теперь вспомним, что Рис. 1 строился для реальных составляющих $S$ и $\alpha$. Но если выражение (6) и функция $S(x)$ аналитические, то обобщение (7) для комплексных $\alpha$, $x_i$ тоже будет полностью удовлетворять условиям экстремума (минимума), если этим условиям удовлетворяет её реальная составляющая. Ведь при этом функции и производные от них сохраняют свои величины, что обеспечивает и сохранение уравнений, образованных из них. Как известно, аналог формулы (7) в римановой геометрии,  опирающийся на принятые представления о параллельном переносе вектора вдоль геодезической линии, имеет вид

\begin{center}
\begin{equation}
    \frac{d^2x_i}{ds^2} + \sum \Gamma_{\alpha \beta}^i (\frac{dx_\alpha}{ds}) (\frac{dx_\beta}{ds}) = 0, 
\end{equation}
\end{center}

\noindent где $\Gamma_{\alpha \beta}^i$ - символ Кристоффеля, довольно громоздко выражающийся через $\rho_{ij}$ и их производные по  $x_v$ . Сейчас (8) стало одной из основ теории векторного и тензорного анализа. Вывод (7), сделанный на другом принципе, гораздо проще, а сама формула прозрачнее, понятнее. 

Какие соображения возникают при сопоставлении этих формул?  Они взаимозаменяемы, но (7) может иметь преимущество перед (8), быть удобнее в использовании. Она может, в частности, служить для иллюстраций и объяснений итогов применения (8). Заметим, что (7) и (8) испытывались на выводе широко известного в радиофизике закона Снелиуса. Обе его подтвердили, причём (7) – сразу, а (8) – после длительных кропотливых вычислений. Очевидно, формула (7) весьма перспективна в смысле применения в теоретической и практической физике. Кроме того, она в силу своей наглядности очень удобна для быстрых качественных оценок перед  стадией окончательных точных расчётов.                                                        Особо отметим следующий факт. Формула (8) в своё время послужила основой для создания общей теории относительности.  Не  исключено, что формула (7) окажется полезной для развития работ и в этом направлении.

\newpage

\begin{figure}[h]
    \centering
    \includegraphics[scale = .75]{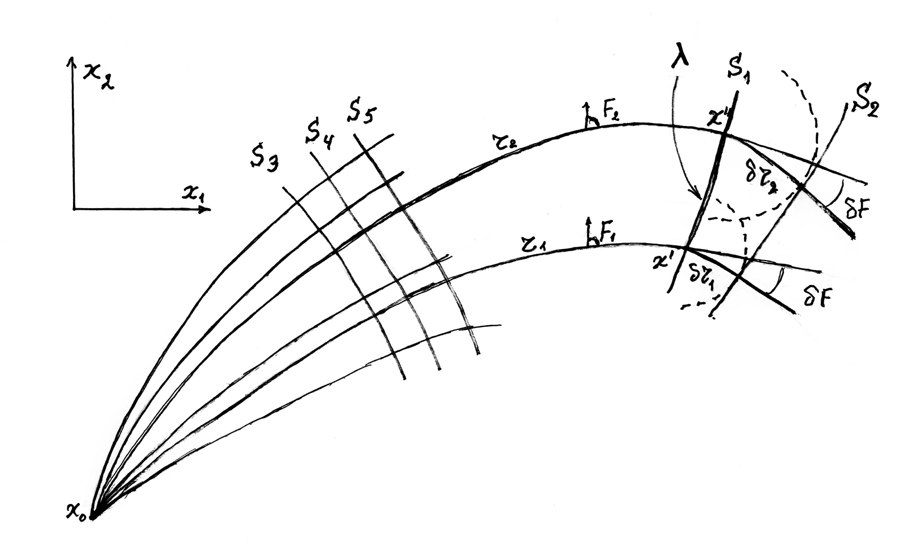}
    \caption{}
\end{figure}

\bibliographystyle{amsplain}